\documentclass[12pt,a4paper]{article}

\usepackage[english]{babel}
\usepackage{graphicx}
\usepackage[centertags]{amsmath}
\usepackage{amsfonts}
\usepackage{amssymb,amsbsy}
\usepackage{amsthm}
\usepackage{newlfont}

\newcommand{\be}{\beta}
\newcommand{\la}{\lambda}

\newcommand{\gf}{\mathfrak G}
\newcommand{\ts}{\tilde S}
\newcommand{\bx}{\mathbf x}

\newcommand{\by}{\mathbf y}
\newcommand{\bz}{\mathbf z}
\newcommand{\bxp}{\mathbf x'}
\newcommand{\bq}{\mathbf q}

\newcommand{\bxo}{\stackrel{o}{\bx}}
\newcommand{\bzo}{\stackrel{o}{\bz}}
\newcommand{\ubx}{\underline\bx}
\newcommand{\obx}{\overline{\bx}}
\newcommand{\ubz}{\underline\bz}
\newcommand{\obz}{\overline{\bz}}

\newcommand{\cx}{\check{x}}

\newcommand{\ey}{\text{\O}}
\newcommand{\bs}{\backslash}
\newcommand{\bcap}{\bigcap}
\newcommand{\bcup}{\bigcup}
\newcommand{\1}{\mathbf 1}
\newcommand{\0}{\mathbf 0}
\newcommand{\lb}{\left(}
\newcommand{\rb}{\right)}

\newcommand{\st}{\subset}

\def\ind#1{{\mbox{\fontsize{9}
{10.8}\selectfont$\relax#1 $} }}

\newtheorem{thm}{Theorem}
\newtheorem{cor}[thm]{Corollary}

\newtheorem{prop}[thm]{Proposition}
\theoremstyle{remark}
\newtheorem{rem}{Remark}

\newcommand{\nn}{\noindent}

\begin{document}

\title{Network flow optimization for restoration of images}

\author{Boris A. Zalesky\\
Institute of Engineering Cybernetics, Minsk, Belarus
\footnote{address: 220012, Surganov str. 6, Minsk,
Belarus}}

\maketitle

\begin{abstract}
The network flow optimization approach is offered for
restoration of  gray-scale and color images
corrupted by noise. The Ising models are used as a
statistical background of the proposed method. The new
multiresolution network flow minimum cut algorithm,
which is especially efficient in identification of the
maximum a posteriori estimates of corrupted images, is
presented. The algorithm is able to compute the MAP
estimates of large size images and can be used in a
concurrent mode. We also consider the problem of
integer minimization of two functions
$U_1(\bx)=\la\sum_i|y_i-x_i|+
\sum_{i,j}\be_{i,j}|x_i-x_j|$
and
$U_2(\bx)=\sum_i\la_i(y_i-x_i)^2+
\sum_{i,j}\be_{i,j}(x_i-x_j)^2$
with parameters
$\la,\la_i,\be_{i,j}>0$
and vectors
$\bx=(x_1,\ldots,x_n), \by=(y_1,\ldots,y_n) \in
\{0,\ldots,L-1\}^n$.
Those functions constitute the
energy ones for the Ising model of color and grayscale
images. In the case $L=2$  they coincide, determining
the energy function of the Ising model of binary
images, and their minimization becomes equivalent to
the network flow minimum cut problem. The efficient
integer minimization of
$U_1(\bx),U_2(\bx)$
by the
network flow algorithms is described.
\end{abstract}

\vspace{2cm}
\nn

\text{}\\
\nn

\newpage

\section{Introduction}

We present a new multiresolution algorithm
for finding the minimum network flow cut  and
methods of efficient integer minimization of the
function
$U_1(\bx)=\la\sum_i |y_i-x_i|+
\sum_{i,j}\be_{i,j}|x_i-x_j|$
and
$U_2(\bx)=\sum_i\la_i(y_i-x_i)^2+
\sum_{i,j}\be_{i,j}(x_i-x_j)^2$.
The problem is posed in terms of Bayesian approach
to image restoration to have unified canvas
of presentation of the results,
and since the  results developed were tested and
turned out efficient for processing
of corrupted images.
Though the minimization of the functions
$U_1(\bx)$ and $U_1(\bx)$ concerns not only
the MAP estimation. It can be understood,
for instance, as $l_1$- and
$l_2$-regression and so on.

The restoration of degraded images is a branch of
image processing that is now extensively studied for
its evident practical importance as well as
theoretical interest. There are many approaches to
solution of the problem.

We consider here the Maximum a posteriori (MAP)
Bayesian estimation \cite{GG84,G90}. As usual, in
this statistical method one supposes an unknown image to be
corrupted by a random noise, and that a corrupted version
of the image is only observable. The unknown original
image is presented as a random realization of
some Markov random field (MRF), which is independent
of the random noise. The distributions of the MRF
and the random noise are assumed to be known. The MAP
estimate is determined as the Maximum likelihood
evaluation of the conditional distribution
of the original image given the observed one.

Because of use of an appropriate prior information,
the Bayesian MAP estimators are among those,
which give the best quality of restored images.
The theory of evaluation
of the MAP estimates is also
significantly developed
\cite{GG84,G90,Gi95,GPS89}.
But one of disadvantages of almost all known
algorithms is their NP-hardness (their execution time is
exponent in a degree of image size).

The first efficient method for evaluation of the
Ising MAP estimate of binary images, which requires
polynomial execution time in number of image pixels,
was proposed by Griege, Porteous and Seheult
\cite{GPS89}.
The authors used
the network flow optimization approach to reduce
problem to identification of the maximum network flow.
The finding of the exact
Ising MAP estimate of binary image required $O(n^3)$
operations, where $n$ was an image size. But known
maximum network flow algorithms still did not allow
computation of the real binary images having size
$n=256\times 256$ or more. Therefore, in \cite{GPS89}
an heuristic modification of the algorithm has been
described.  That algorithm was implemented not for
entire images but for  pieces of images with fixed
values at their frontiers.

The similar idea was used for the proposed
\emph{multiresolution network flow minimum cut
algorithm}
(actually, to find the MAP of a corrupted image  one  need
only the minimum network flow cut).
Some additional theoretical reasons allowed describing
a new algorithm, which identifies the  exact minimum cut using
special modifications of subnetworks of the
partitioned original network. This algorithm can be
exploited for an arbitrary network. At worst
it takes $O(n^3)$ operations for computation of the
minimum network flow cut, but for networks that
correspond to real corrupted images it usually requires only
$O(n)$ operations, which would be implemented in a concurrent
mode.

The description, theoretical ground
and the use of the multiresolution network
flow minimum cut algorithm (MNFC) is presented in
Section~\ref{s:MNFC}.

The brief entry to two Ising models of color and
grayscale images is done in Section~\ref{s:Ising}.
Those models are widely used for reconstruction of
corrupted images, they coincide when images are
binary \cite{GG84,G90,Gi95}.
In the same Section
the problem of finding the exact MAP estimates
is formulated in terms of the integer programming.

In Section~\ref{s:min} methods of the integer
minimization of functions (which are the energy
functionals of the considered Ising models)
$U_1(\bx)=\la\sum_i|y_i-x_i|+
\sum_{i,j}\be_{i,j}|x_i-x_j|$
and
$U_2(\bx)=\sum_i\la_i(y_i-x_i)^2+
\sum_{i,j}\be_{i,j}(x_i-x_j)^2$
with parameters
$\la,\be_{i,j}>0$
and vectors
$\bx=(x_1,\ldots,x_n), \by=(y_1,\ldots,y_n) \in
\{0,\ldots,L-1\}^n$ by the network flow algorithms
(including the MNFC) is described.

Applications of developed methods is
given in Section~\ref{s:appl}.

\section{Multiresolution network flow minimum cut
algorithm}\label{s:MNFC}

\subsection
{ Preliminaries and Notations \label{s:P&N}}
We now describe the  \emph{multiresolution network flow
minimum cut algorithm}. For some types of
networks this highly parallel method turns out more
speedy and more efficient in comparison with known
maximum network flow algorithms. It was successfully
used for identification of  minimum cuts of  large networks
(for instance, for determination of the MAP of binary
images \cite{GPS89})
while classical methods were not able to solve the
problem.

The essence of the MNFC is the following. A network
$\gf$ is partitioned into several subnetworks of
appropriate sizes. Every subnetwork is modified in a
special way and evaluated two times. The first time
all boundary arcs going from the outside into
the subnetwork are considered as going from the source
(to the same nodes) and all boundary arcs going
from within the subnetwork outside are ruled out the
consideration. The minimum cut of the modified
subnetwork is determined by known algorithm.
It is proven that all nodes, which are connected
with the sink by a directed path,
will be connected with the sink in
the solution of the original network by the same path.
The second time all boundary arcs going from
the outside into the subnetwork are excluded from the
consideration and all boundary arcs going
from within the subnetwork outside are supposed to be
connected with the sink. The minimum cut of the
modified subnetwork is identified once more. This time all
directed paths going from the source are determined.
The paths identified are excluded from the further
consideration.
The reduced network $\gf_1$ is divided into several
subnetworks again. The procedure of identification of
arcs belonging to  $\gf_1$ but not belonging to the
minimal cut  is repeated now with one difference --
we take in account arcs connected with excluded nodes.
The arcs not belonging to the minimal cut of $\gf_1$
are found and are removed anew. We go to higher levels
until we obtain the network $\gf_k$ that can be solved
by a usual maximum network flow algorithm.

In more details. Let the network $\gf$
consist of $n+2$ numbered \emph{nodes}
$\ts=\{0,1,2\ldots,n,n+1\}$, where $s=0$ is the
\emph{source}, $t=n+1$ is the \emph{sink} and $
S=\{1,2\ldots,n\}$ are usual nodes. The set of
directed \emph{arcs} is $A=\{(i,j)\ :\ i,j\in \ts\}$.
Capacities of arcs are denoted by $d_{i,j}>0$.
Suppose the network $\gf$ satisfies the condition:

\nn (G) \hspace{0.1cm}\hfill \parbox{12.5cm}{ For every
usual node $i\in S$ there is either an arc $(s,i)\in
A$ connecting the source $s$ with $i$ or  an arc
$(i,t)\in A$ connecting $i$ with the sink $t$ but not both
of these arcs.}

\begin{rem}
The condition (G) does not restrict the use of the MNFC. It
has been done to simplify writing down notations and
proofs. Any network $\gf$ can be easily
modified to satisfy (G) for  $O(n)$ operations (one
look through usual network nodes $S$). The modified
network will have the same minimum cut.
\end{rem}

Let two sets $W,B\subset S$ such that $W\bcap B=\ey$
and $W\bcup B=S$ be a partition of $S$.
The vector $\bx\ : S\rightarrow
\{0,1\}^{|S|}$ with coordinates $x_i=1$, if $i\in W$,
and  $x_i=0$ otherwise is the indicator vector of
the set $W$. The set of all such indicator vectors $\bx$
is denoted by $\mathbf X$. The \emph{capacity of the
$s-t$ cut} $C(\bx)$
is defined as the sum of capacities of all
forward arcs going from the set $W\cup \{s\}$ to the
set $B\cup \{t\}$ \cite{PapSteig82}, i.e.
$$C(\bx)=\sum\limits_{i\in W\cup \{s\}}
\sum\limits_{j\in B\cup \{t\}}d_{i,j}.$$
Let $\by$ be the vector with coordinates $y_i=1$, if
$(s,i)\in A$, or $y_i=0$, if $(i,t)\in A$ (remind that
the condition (G) is valid). Set
$\la_{i}=d_{s,i}$, if $(s,i)\in A$, or
$\la_{i}=d_{i,t}$, if $(i,t)\in A$. The capacities of
usual arcs $(i,j),\ i,j\in S$ will be denoted by
$\be_{i,j}=d_{i,j}$.
Then in new notations
$$C(\bx)=\sum_{
\stackrel{\mbox{\fontsize{8}{9.6}\selectfont$ i\in
W,$}}
{\mbox{\fontsize{8}{9.6}\selectfont${y_i=0}$}}}\la_i+
\sum_{ \stackrel{\mbox{\fontsize{8}{9.6}\selectfont$
i\in B,$}}
{\mbox{\fontsize{8}{9.6}\selectfont$y_i=1$}}}\la_i\ +
\ \sum_{(i,j)\in S\times S}\be_{i,j}(x_i-x_j)x_i.$$
Now introduce the function
$$ U(\bx)=\sum_{i\in
S}\la_i (1-2y_i)x_i\ + \ \sum_{(i,j)\in S\times
S}\be_{i,j}(x_i-x_j)x_i. $$
It is easy to see that
$C(\bx)=U(\bx)+\sum_{i=1}^n\la_i y_i$. Since the term
$\sum_{i=1}^n\la_iy_i$ does not depend on $\bx$, the
functions $C(\bx)$ and $U(\bx)$ has minima in the same
points $\bx^*=\text{argmin}_{\bx\in X} C(\bx)=
\text{argmin}_{\bx\in X} U(\bx)$. Therefore,
the solutions $\bx^*=\text{argmin}_{\bx\in X} U(\bx)$
entirely identify minimum network flow cuts
\cite{FF62,PR75}.

For an arbitrary subset of usual nodes $E\subset S$
we define two functions
\begin{eqnarray}\label{e:U}
{ U_{E}(\bx)=\sum_{i\in E}\la_i (1-2y_i)x_i\
{}}
+\sum_{(i,j)\in (E\times E)\cup (E\times E^c) \cup (E^c
\times E)} \be_{i,j}(x_i-x_j)x_i
\end{eqnarray}
and
$$V_{E}(\bx)=\sum\limits_{i\in E^c} \la_i (1-2y_i)x_i\
+\ \sum\limits_{(i,j)\in (E^c \times E^c)}
\be_{i,j}(x_i-x_j)x_i,$$
the sum of which is equal to
$U(\bx)=U_{E}(\bx)+V_{E}(\bx)$. Also we define
restriction $\bx_{E}$  of $\bx$ onto the set $E$ such that
$\bx=(\bx_{E},\bx_{E^c})$. The function $V_{E}(\bx)$
depends only on $\bx_{E}$ (i.e.
$V_{E}(\bx)=V_{E}(\bx_E)$), therefore the following simple
proposition is valid.
\begin{prop}\label{p:locmin}
If  $\bxp_{E}$ minimizes the function
$\phi(\bx_{E})=U(\bx_{E},\bxo_{E^c})$ for some fixed
$\bxo_{E^c}$, then for any set  $D\st E$ the
restriction  $\bxp_{D}$ of $\bxp_{E}$ minimizes the
function $\psi(\bx_{D})=
U_{D}(\bx_{D},\bxp_{E\backslash D},\bxo_{E^c})$ and
vice versa.
\end{prop}

\nn For vectors  $\bx_E,\bz_E$ we will write:

$\bx_E\le \bz_E$, if  $x_i\le z_i,\ i\in E$;

$\bx_E<\bz_E$, if $x_i\le z_i,\ i\in E$ and there is
at least one node  $j\in E$ such that  $x_j<z_j$;

$\bx_E\nleq\bz_E$, if there are nodes  $i,j\in E$ such
that  $x_i < z_i$ and  $x_j > z_j$.

\nn Our method is based on
monotone (in some sence) dependence of restrictions
$\bx^*_E$ of solutions
$\bx^*=\text{argmin}_{\bx\in X} U(\bx)$ on
values $\bx^*_{E^c}$.

\subsection{Properties of Minimum Cuts\label{s:PMC}}
Let us consider the property of the  monotony
of $\bx^*_E$ in details.
For fixed  $\bxo_{E^c}$ and $\bzo_{E^c}$ set
$\bx'_E=\text{argmin}_{\bx_E}U(\bx_{E},\bxo_{E^c})$
and
$\bz'_E=
\text{argmin}_{\bx_E}U(\bx_{E},\bzo_{E^c})$.
It is clear that every set $\{\bx'_E\}_{\bxo_{E^c}}$,
$\{\bz'_E\}_{\bzo_{E^c}}$ can consist of more than one
element.
In general,
in the case  $\bxo_{E^c}\le\bzo_{E^c}$
the inequality
$\bx'_E\le\bz'_E$ is not satisfied,
but without fail
there exist  solutions $\bx'_{E}$ and  $\bz'_{E}$ such
that $\bx'_{E}\le\bz'_E$.
\begin{thm}\label{t:monot}
Suppose that
${\stackrel{o}{\bx}}_{E^c}\le\stackrel{o}{\bz}_{E^c}$
are fixed and let $\bx'_E\nleq\bz'_E$ (i.e. solutions
$\bx'_E$ and $\bz'_E$ are not comparable). Then for
nonempty set $D=\{i\in E\ :\ x'_i=1,\ z'_i=0\}$ the
vector $(\bz'_D,\bx'_{E\bs D})=(\0_D,\bx'_{E\bs D})$
minimizes the
function $U$ under condition  $\bxo_{E^c}$, and the
vector $(\bx'_D,\bz'_{E\bs D})=(\1_D,\bz'_{E\bs D})$
minimizes the
function $U$ under condition  $\bzo_{E^c}$, that is
$(\bz'_D,\bx'_{E\bs D})\in\{\bx'_E\}_{\bxo_{E^c}}$,
$(\bx'_D,\bz'_{E\bs D})\in\{\bz'_E\}_{\bzo_{E^c}}$ and
$(\bz'_D,\bx'_{E\bs D})\le (\bx'_D,\bz'_{E\bs D})$.
\end{thm}
\begin{proof}
Suppose $\bx'_E\in\{\bx'_E\}_{\bxo_{E^c}}$ and
$\bz'_E\in\{\bz'_E\}_{\bzo_{E^c}}$ are two
incomparable solutions  $\bx'_E\nleq\bz'_E$ for
frontier conditions $\bxo_{E^c}\le\bzo_{E^c}$.
It follows from Proposition \ref{p:locmin}, equalities
$\bx'_D=\mathbf 1_D$ and $\bz'_D=\mathbf 0_D$, and the
inequality $(\bx'_{E\bs D},\bxo_{E^c})\le (\bz'_{E\bs
D},\bzo_{E^c})$ that
\begin{equation}\label{e:in1}
U_D(\bz'_D,\bz'_{E\bs D},\bzo_{E^c})\le
U_D(\bx'_D,\bz'_{E\bs D},\bzo_{E^c})
\le U_D(\bx'_D,\bx'_{E\bs D},\bxo_{E^c}).
\end{equation}
And by analogy,
\begin{equation}\label{e:in2}
U_D(\bx'_D,\bx'_{E\bs D},\bxo_{E^c})\le
U_D(\bz'_D,\bx'_{E\bs D},\bxo_{E^c})
\le U_D(\bz'_D,\bz'_{E\bs D},\bzo_{E^c})
\end{equation}
Gathering  estimates (\ref{e:in1},\ref{e:in2}) in one
chain we can see that all inequalities in
(\ref{e:in1},\ref{e:in2}) are actually equalities.
This fact and Proposition \ref{p:locmin} finish the
proof.
\end{proof}
The following Corollary, which
characterizes some structural properties of the set of
solutions $\{\bx^*\}$, is deduced from Theorem
\ref{t:monot}.
\begin{cor}\label{c:struct}
(i) If fixed frontier vectors satisfy the condition
${\stackrel{o}{\bx}}_{E^c}\le\stackrel{o}{\bz}_{E^c}$,
then for any solution
$\bx'_E=\text{argmin}_{\bx_E}U(\bx_{E},\bxo_{E^c})$
there exists a solution
$\bz'_E=\text{argmin}_{\bx_E}U(\bx_{E},\bzo_{E^c})$
such that $\bx'_E\le\bz'_E$.

(ii) If fixed frontier functions  satisfy the
condition
${\stackrel{o}{\bx}}_{E^c}\le\stackrel{o}{\bz}_{E^c}$,
then for any solution
$\bz'_E=\text{argmin}_{\bx_E}U(\bx_{E},\bzo_{E^c})$
there exists a solution\\
$\bx'_E=\text{argmin}_{\bx_E}U(\bx_{E},\bxo_{E^c})$
such that $\bx'_E\le\bz'_E$.

(iii) For any frontier condition
${\stackrel{o}{\bx}}_{E^c}$ the set
$\{\bx'_E\}_{\bxo_{E^c}}$ has  the minimal (the maximal)
element  $\ubx'_{E}$ ($\obx'_{E}$).

(iv) If
${\stackrel{o}{\bx}}_{E^c}\le\stackrel{o}{\bz}_{E^c}$,
then $\ubx'_{E}\le\ubz'_{E}$ and
$\obx'_{E}\le\obz'_{E}$.
\end{cor}

Sentences (\emph{i}) and (\emph{ii}) follow
immediately
from Theorem~\ref{t:monot}. Sentence (\emph{iii})
is deduced from Theorem~\ref{t:monot} for
$\bxo_{E^c}=\bzo_{E^c}$. And (\emph{iv}) follows from
(\emph{i,ii}) and definitions of  minimal and the maximal
elements.

To identify a local solution
$\bx'_E=\text{argmin}_{\bx_E}U(\bx_{E},\bxo_{E^c})$
 the network  $\gf'_E$ with nodes $\tilde
S_E=\big\{E,\{s\},\{t\}\big\}$ and arcs $(i,j),\
i,j\in\tilde S$ having capacities
\begin{equation}\label{e:capac}
d_{i,j}=\be_{i,j},\ \text{если} \ i,j\in E;\quad
d_{s,i}=\la_{i}y_i+ \sum\limits_{\stackrel{j\in
E^c}{x_j=1}}\be_{j,i};\quad
d_{i,t}=\la_{i}(1-y_i)+
\sum\limits_{\stackrel{j\in E^c}{x_j=0}}\be_{i,j}
\end{equation}
is used, since the following proposition can be easily
proved.
\begin{prop}\label{p:locsol}
The local solution $\bx'_E$ sets the minimum cut of the
network $\gf'_E$ and vice versa.
\end{prop}
\nn
The subsets $E\in S$ are chosen so that the maximum flow
in the network $\gf'_E$ can be computed by  usual
maximum network flow algorithms.

\subsection{The MNFC Algorithm\label{s:A}}
The main idea of the MNFC is
to estimate at least  parts of restrictions
$\bx^*_{E_\ind{i}}$ of  solutions
$\bx^*=\text{argmin}_\bx U(\bx)$ for a suitable
partition $\bcap E_\ind{i}=S,\ E_\ind{i}\bcap E_\ind{j}=\ey$.
For this purpose the monotone
dependence of  local solutions $\bx^*_{E_\ind{i}}$
on the frontier
values of $\bx^*_{E^c_\ind{i}}$
(in the sense of Corollary \ref{c:struct})
is exploited. The parts of $\bx^*_{E_\ind{i}}$ estimated by
the special local solutions $\bx^{'}_{E_\ind{i}}$
are ruled out the further consideration. It
significantly reduces computational expenses.

More precisely. Let sets
$E_1(1),\ldots,E_{k_\ind{1}}(1)$ partition the
set $S$ so that $\bcap^{k_\ind{1}}_{i=1}E_i(1)=\ey$
and $\bcup^{k_\ind{1}}_{i=1} E_i(1)=S$. It follows
from Proposition \ref{p:locmin} that
$\bx'_{E_\ind{i}(1)}=\text{argmin}_{\bx_{E_\ind{i}(1)}}
U(\bx_{E_\ind{i}(1)},\bx^*_{E^c_\ind{i}(1)})$
$=\bx^*_{E_{\ind{i}}(1)}$, and the solution
$\bx'_{E_\ind{i}(1)}$ minimizes the function
$U_{E_\ind{i}(1)}
(\bx_{E_\ind{i}(1)},\bx^*_{E^c_\ind{i}(1)})$ for the
frontier condition  $\bx^*_{E^c_\ind{i}(1)}$.
Corollary \ref{c:struct} guarantees existence of
two solutions (that can be found by known
maximum network flow algorithms)
$$\bx_{0,E_\ind{i}(1)}=
\text{argmin}_{\bx_{E_\ind{i}(1)}}
U(\bx_{E_\ind{i}(1)},\0_{E^c_\ind{i}(1)})$$
and
$$\bx_{1,E_\ind{i}(1)}=
\text{argmin}_{\bx_{E_\ind{i}(1)}}
U(\bx_{E_\ind{i}(1)},\1_{E^c_\ind{i}(1)})$$
satisfying the inequality
\begin{equation}\label{e:x01}
\hspace{-0.7cm}\bx_{0,E_\ind{i}(1)}\le
\bx'_{E_\ind{i}(1)}\le\bx_{1,E_\ind{i}(1)}
\end{equation}
Denote the sets on nodes
$B_\ind{i}(1)=\{k\in E_\ind{i}(1)\ |\
x_{1,E_\ind{i}(1),k}=0\}$ and
$W_\ind{i}(1)=\{k\in E_\ind{i}(1)\ |\
x_{0,E_\ind{i}(1),k}=1\}$.
The  equalities
$\0_{B_\ind{i}(1)}=\bx_{1,B_\ind{i}(1)}=
\bx^*_{B_\ind{i}(1)}$ and
$\1_{W_\ind{i}(1)}=\bx_{0,W_\ind{i}(1)}=
\bx^*_{W_\ind{i}(1)}$
is inferred from (\ref{e:x01}).

If the set
$R(1)= \bcup^{k_\ind{1}}_{i=1} (W_\ind{i}(1)\bcup
B_\ind{i}(1))$ is not empty, then we identified the
part of the solution  $\bx^*_{R(1)}$. There is a sense to
continue  the MNFC . Assign  $S(2)=S\bs R(1)$ and
consider now
the reduced problem of identification
$\bx^*_{S(2)}=\text{argmin}_{\bx_{S(2)}}
U(\bx_{S(2)},\bx^*_{R(1)})=$
$\text{argmin}_{\bx_{S(2)}}
U_{S(2)}(\bx_{S(2)},\bx^*_{R(1)})$
for the frontier condition $\bx^*_{R(1)}$.
Partition
$S(2)=\bcup^{k_\ind{2}}_{i=1} E_i(2),
E_i(2)\bcap E_j(2)=\ey$ and estimate
$$\bx'_{E_\ind{i}(2)}=\bx^*_{E_{\ind{i}}(2)}=
\text{argmin}_{\bx_{E_\ind{i}(2)}}
U_{E_\ind{i}(2)}(\bx_{E_\ind{i}(2)},\bx^*_{S(2)\bs
E_\ind{i}(2)},\bx^*_{R(1)})$$
by the solutions
$$\bx_{1,E_\ind{i}(2)}=
\text{argmin}_{\bx_{E_\ind{i}(2)}}
U_{E_\ind{i}(2)}(\bx_{E_\ind{i}(2)},\1_{S(2)\bs
E_\ind{i}(2)},\bx^*_{R(1)})$$
and
$$\bx_{0,E_\ind{i}(2)}=
\text{argmin}_{\bx_{E_\ind{i}(2)}}
U_{E_\ind{i}(2)}(\bx_{E_\ind{i}(2)}
\0_{S(2)\bs E_\ind{i}(2)},
\bx^*_{R(1)})$$
satisfying the inequality
$\bx_{0,E_\ind{i}(2)}\le\bx'_{E_\ind{i}(2)}
\le\bx_{1,E_\ind{i}(2)}$.
Then, consider sets
$B_\ind{i}(2)=\{k\in E_\ind{i}(2)\ |\
x_{1,E_\ind{i}(2),k}=0\}$
and $W_\ind{i}(2)=\{k\in E_\ind{i}(2)\ |\
x_{0,E_\ind{i}(2),k}=1\}$,
which identify
$\bx^*_{B_\ind{i}(2)}=\0_{B_\ind{i}(2)}$ and
$\bx^*_{W_\ind{i}(2)}=\1_{W_\ind{i}(2)}$ correspondingly.
If the set
$R(2)=\bcup^{k_\ind{2}}_{i=1}
(W_\ind{i}(2))\bcup B_\ind{i}(2))$
is nonempty, the algorithm is employed at  the
third level.
The MNFC is iterated at higher levels
until the problem will be completely solved  or
until $R(l)\neq\ey$.
At the uppermost level
an appropriate known  \emph{maximum network flow
algorithm}
(MNF) is applied.

Let us describe the MNFC in brief.

\nn
{\bf Step 1} (Initialization).
Assign: the level number  $l=1$, the set not
estimated nodes  $S(l)=S$, the set of estimated nodes
$R=\ey$, and the set of arcs $A(l)=A$.

\nn
{\bf Step 2} (Partition of the network).
If the maximum flow in the
network $\{S(l),\{s\},\{t\},$ $A(l)\}$ can  be efficiently
identified by usual MNF, we do not need its partition.
Assign the number of partition elements  $k_l=1$, and
$\ E_1(l)=S(l)$. Identify the maximum flow in
$\{S(l),\{s\},\{t\},A(l)\}$,  \emph{STOP}.
Otherwise, partition the set $S(l)$ by not intersected
sets
$E_1(l),E_2(l),\ldots,$ $E_{k_\ind{l}}(l),$
$\bcap^{k_\ind{l}}_{i=1}E_i(l)=\ey,\
\bcup^{k_\ind{l}}_{i=1} E_i(l)=S(l)$,
so that every network
$\{E_i(l),$ $\{s\},\{t\},A_i(l) \}$ with arcs
$A_i(l)=A_{E_\ind{i}}(l)=\{(\mu,\nu)\ :\ \mu\in E_i(l)\
\text{or}\ \nu\in E_i(l) \}$
can be evaluated by an appropriate maximum network flow
algorithm (for instance,
by the labeling algorithm or other ones).

\nn
{\bf Step 3} (Computation of local estimates
$\bx_{0,E_\ind{i}(l)}, \bx_{1,E_\ind{i}(l)}$).
Compute the vectors
\begin{multline*}
\bx_{0,E_\ind{i}(l)}=
\text{argmin}_{\bx_{E_\ind{i}(l)}}
U(\bx_{E_\ind{i}(l)},\0_{S(l)\bs E_\ind{i}(1)},
\bx^*_R)\\
\le
\bx_{1,E_\ind{i}(l)}=
\text{argmin}_{\bx_{E_\ind{i}(l)}}
U(\bx_{E_\ind{i}(l)},\1_{S(l)\bs E_\ind{i}(l)},
\bx^*_R)
\end{multline*}
(see Corollary \ref{c:struct} and Proposition
\ref{p:locsol})
by an appropriate MNF.
These two vectors correspond to minimum  cut of the
networks $\{E_i(l),\{s\},\{t\},A_i(l)\}$ for
frontier conditions
$(\0_{S(l)\bs E_\ind{i}(1)},\bx^*_R)$ and
$(\1_{S(l)\bs E_\ind{i}(1)},\bx^*_R)$ respectively.

\nn
{\bf Step 4} (Identification of the set
$R(l)$ of estimated nodes).
Find the sets of nodes
$B_\ind{i}(l)=\{i\in E_\ind{i}(l)\ | \
x_{1,E_\ind{i}(l),i}=0\}$
and
$W_\ind{i}(l)=\{i\in E_\ind{i}(l)\ | \
x_{0,E_\ind{i}(l),i}=1 \}$
keeping their values in  $\bx^*$. Assign the set
$R(l)=\bcup^{k_\ind{l}}_{i=1}
(W_\ind{i}(l)\bcup B_\ind{i}(l))$

\nn
{\bf Step 5} (Check whether  the multiresolution
approach can be continued).
If   $R(l)=\ey$, interrupt execution of
the MNFC and try to identify $\bx_{S(l)}^{*}=$
$\text{argmin}_{\bx_{S(l)}}U(\bx_{S(l)},\bx^*_R )$
by an appropriate MNF,
\emph{STOP}.

\nn
{\bf Step 6} (Jump to the higher level).
Assign  $R:=R\bigcup R(l)$ and
$S(l+1)=S(l)\bs R$. If  $S(l+1)=\ey$, the proplem is
solved -- \emph{STOP}. Otherwise, go to the higher
level, i.e. assign  $l:=l+1$,
$A(l)=\{(\mu,\nu)\ :\ \mu\in S(l)\
\text{or}\ \nu\in S(l) \}$
and go to {\bf  Step 2.}

The {\bf Step 3}  of the MMC can be efficiently
executed in a concurrent mode.

\subsection{The number of operations and execution time}

The number of operations required to execute the
MNFC essentially depends on properties of the initial
network $\gf$. In the case $R(1)=\ey$ the MNFC is
reduced to the usual maximum network flow (MNF) algorithm.
But there are a lot of
networks admitting efficient application of the MNFC.
Applications of the algorithm
to the decision of  real problems are described
in Section \ref{s:appl}.
In those applications identification of  network minimum
cuts have required only $O(n)$ operations.

In the following  Proposition we formulate
the simple necessary and sufficient condition in order
to the MNFC does not degenerate to the usual maximum
network flow algorithm.
\begin{prop}
The set $W_\ind{i}(l)\neq \ey$,
respectively, $B_\ind{i}(l)\neq \ey$
if and only
if there exists such a set $D\subset E_\ind{i}(l)$ for
which the inequality
\begin{multline}\label{e:necsaf}
\sum_i\la_i(1-2y_i)+\sum_{(i,j)\in (D\times
D^c)}\be_{i,j}\le 0,\quad\text{respectively},\\
\sum_i\la_i(1-2y_i)+\sum_{(i,j)\in (D^c\times
D)}\be_{i,j}\ge 0
\end{multline}
is valid.
\end{prop}

We will write down  estimates of operations number
$N_{\text{op}}$ and  time of parallel execution
$T_{\text{par}}$ in
general form and consider several special cases.
Let $n_i(l)$ be number of nodes in a subnetwork,
$m_i(l)$ be number of arcs in a modified subnetwork,
and $a_i(l)$ be number
of boundary arcs such that
$a_i(l)=(\mu,\nu),\ \mu\in E_i(l),\nu\in E^c_i(l)\
\text{or}\ \nu\in E_i(l),\mu\in E^c_i(l)$.
It is supposed that the network satisfies the condition
(G) (recall that modification of the network
to satisfy this condition needs $O(n)$ operations).
The amount of operations $N_{\text{op}}$ has the same order
as  number of operations required to
construct all local estimates
$\bx_{0,E_\ind{i}(l)},\bx_{1,E_\ind{i}(l)}$.
But the computation of each local estimate requires
$O(a_i(l))+O(n^3_i(l))$ operations, where
$O(a_i(l))$ operations are needed to modify  subnetwork
capacities (see, (\ref{e:capac})) and $O(n^3_i(l))$
operations are required directly to identify the local
estimates (the second term  $O(n_i(l)^3)$ depends on the used
maximum network flow algorithm. It can be replaced by
$O(n_i(l)m^2_i(l))$ etc.). If the number of levels
required to compute the minimum cut is equal to $L$,
then the number of operations is
$N_{\text{op}}=O\lb\sum^L_{l=1}
\sum_{i=1}^{k_l}(a_i(l)+n^3_i(l))\rb$
and parallel execution time is
$T_{\text{par}}=O\lb\sum^L_{l=1}
\max_{i}(a_i(l)+n^3_i(l))\rb $.

First we consider the worst case. Suppose that all
usual nodes $S$ are completely connected, i.e. every
two usual nodes are connected by two directed arcs.
And suppose that a prior information about a strategy
of partitions is not known.
Then, the pyramidal approach allows to
avoid excessive computational expenses. At level
$1\le l\le \log_2 n=L$
partition the original network into
$2^{k(l)},\ k(l)=\log_2 n-l$
subnetworks containing the same number of nodes.
Therefore, at worst, the number of operations is
$N_{\text{op}}=O( n^3)$ and the parallel
execution time is $T_{\text{par}}=O(n^3)$.
The amount of operations will be  $O( n^3)$ if
at each level the small number of nodes is
identified only, and the others are
determined at the uppermost level.
Even in this case computational expenses
of the  MNFC are not excessive in comparison with
the traditional MNF algorithms.

The situation changes if there are  large subnetworks
that satisfy the condition (\ref{e:necsaf}). The
MNFC becomes very efficient. The number
of operation can amount to $O(n)$.
Let us, for instance, consider networks that arise
from the Ising MAP estimation.
Suppose that usual nodes
$S$ form the finite 2D lattice $S=\{1,\cdots,N\}\times
\{1,\cdots,N\}$, and that  every usual node has 2D
index $i=(i_1,i_2)$. The nearest-neighbor nodes are
connected by two differently directed arcs $(i,j)$ and
$(j,i)$ with the same capacity $\be_{i,j}=\be$. Beside
that suppose that each node are connected either with
the source $s$ or with the sink $t$ by directed arcs having
capacity $\la$. Obviously, in this case the condition
(\ref{e:necsaf}) will be satisfied for subnetworks
with nodes having constant values in $\by$
and such that the set of these nodes possesses  small enough
isoperimetric ratio.  Networks that
correspond to corrupted binary images usually contain
large (or even, so called, gigantic) submetworks
satisfying mentioned above property. Those
subnetworks are often estimated by the MNFC and
ruled out the further consideration at the first level of
execution of the algorithm.
So that evaluation of the MAP estimate
takes fractions of a seconds while the classical
maximum network flow algorithms can not identify the
solution at all.
To justify this arguments we formulate the following
sentence.
\begin{prop} \mbox{\fontshape{n}\selectfont (i)}
Let $B_{\by}=\{i\in S\ | \ y_i=0\}$ or
$W_{\by}=\{i\in S\ | \ y_i=1\}$ be connected component
of usual nodes of  the constant value in $\by$.
If $\be\le \la/2$ and $B_{\by}$ (resp., $W_{\by}$) has
a contour,  $B_{\by}$ (resp.,$W_{\by}$)  maintains
the same value in the minimum cut $\bx^*$.

\nn
\mbox{\fontshape{n}\selectfont (ii)}
If for some $M>0$ the capacities satisfy the
inequality $\be\le\frac{M}{2(M+1)}\la$ then every
connected component $B_{\by}$ (resp., $W_{\by}$)
of constant values in $\by$,
consisting of more than $M$ nodes, maintains
the same value in the minimum cut $\bx^*$.
\end{prop}
In the case  $\be > \la/2$
the MNFC also easily finds the MAP estimates of
real binary images (that is equivalent to identification
the minimum cut of the large network).
The results of the practical restoration of corrupted
images by the MNFC are placed in Section \ref{s:appl}.

\section{Two Ising models of color
         and grayscale images}\label{s:Ising}

We describe two simplest Ising models
of color and grayscale images with  energy
functions
$U_1(\bx)=\la\sum_i|y_i-x_i|+
\sum_{i,j}\be_{i,j}|x_i-x_j|$
and
$U_2(\bx)=\sum_i\la_i(y_i-x_i)^2+
\sum_{i,j}\be_{i,j}(x_i-x_j)^2$
with parameters
$\la,\la_i,\be_{i,j}>0$. Those models are well
investigated
\cite{GG84,G90,Gi95,FGN01}
and quite often
applied for practical image restoration.
The main our goal is the developing
efficient algorithms of their integer
minimization.

Let now $S=\{1,\ldots,n\}\times\{1,\ldots,n\}$ denote the finite
square lattice of real numbers; then $i=(i_1,i_2)\in S$  denotes
points or pixels of our images and later on nodes of appropriate
networks. The \emph{gray scale image} is a matrix $\bx_{gray}$
with nonnegative integer coordinates $x_i\in \{0,\ldots,L-1\},
L>2$. The term \emph{"binary image"} is used if coordinates of
the matrix $\bx_{bin}$ take two values $x_i\in  \{0,1\}$.
The \emph{color images} $\bx_{col}$ are specified by three
matrices $(\bx_{r},\bx_{g},\bx_{b})$.

First, we remind in brief of the Ising model of
corrupted binary images
\cite{GG84,G90,Gi95,FGN01}.
This Bayesian approach specifies
an \emph{a priori distribution}
$$
p(\bx)=
c\cdot\exp\big\{-\sum_{i,j\in S}\be_{i,j}(x_i-x_j)^2\big\}
$$
(with the normalizing constant  $c$ and
$\be_{i,j}\ge 0$) over all binary images
$\bx=\bx_{bin}$.
The conditional probability of the original binary image
$\bx_{bin}$ given the corrupted version
$\by=\by_{bin}$ is represented in the form
\begin{equation}\label{e:p(y|x)}
p(\bx|\by)=d(\by)\cdot\exp\big\{-\sum_{i\in
S}\la_i(y_i-x_i)^2-
\sum_{i,j\in S}\be_{i,j}(x_i-x_j)^2 \big\},
\end{equation}
where the function $d(\by)$ is independent of $\bx$.
The MAP estimate, which is defined as a mode of the conditional
distribution $\bx^*=\text{argmax}_\bx p(\bx|\by)$, is
equal to
\begin{equation}\label{e:mode}
\bx^*=\text{argmin}_
{\fontsize{10}{12}\selectfont{\bx}}
\bigg\{ \sum_{i\in S}\la_i(y_i-x_i)^2+
\sum_{i,j\in S}\be_{i,j}(x_i-x_j)^2\bigg\}
\end{equation}

\begin{rem}
Often, to recover real binary images it is supposed
$\la_i=\la$ and $\be_{i,j}>0$ only for $i,j$ belonging to
local neighborhoods (under the last condition the a prior
distribution $p(\bx)$ becomes the Markov Random Field).
We knowingly consider more
complicated model, since even in this case there is an
efficient solution of the problem.
\end{rem}
For binary images the functions $U_1(\bx)$ and $U_2(\bx)$
coincide
and the sum in
(\ref{e:p(y|x)},\ref{e:mode}) is equal to both
of them. But in the case of grayscale images
$U_1(\bx)\neq U_2(\bx)$. Therefore,
for $\bx=\bx_{gray},\by=\by_{gray}$
we consider two
different Bayesian models with the a posteriori
probabilities
\begin{equation}\label{e:p1}
p_1(\bx|\by)=d_1(\by)\cdot\exp\big\{-U_1(\bx) \big\}
\end{equation}
and
\begin{equation}\label{e:p2}
p_2(\bx|\by)=d_2(\by)\cdot\exp\big\{-U_2(\bx) \big\}.
\end{equation}
Both of these models are interesting from
the theoretical and applied points of view. The first one
is less investigated theoretically but it gives
better quality of restored images. Identification of
the MAP for the first model is equivalent to finding
$\bx_1^*=\text{argmin}_
{\fontsize{10}{12}\selectfont{\bx}}
U_1(\bx)$, and for the second --
$\bx_2^*=\text{argmin}_
{\fontsize{10}{12}\selectfont{\bx}}
U_2(\bx)$, where {\it argmin} is taken over all
grayscale images.

For color images
$\bx_{col}=(\bx_{r},\bx_{g},\bx_{b})$
we consider the Ising models
(\ref{e:p1},\ref{e:p2})
for each component separately
(with
$\bx=\bx_{r},\bx=\bx_{g},\bx=\bx_{b}$
and respectively
$\by=\by_{r},\by=\by_{g},\by=\by_{b}$).
Thus, the problem is the same: to find $\bx_1^*$
and $\bx_2^*$ for every component of a color image.

\section{Integer minimization of $U_1(\bx)$ and
 $U_2(\bx)$ by network flow algorithms}\label{s:min}

The network flow methods can be successfully applied
for the efficient integer minimization
of $U_1(\bx)$ and $U_2(\bx)$.
In particular, the MNFC, which has been described in
Section~\ref{s:MNFC},
turns out very efficient for finding of
integer solutions
$\bx_1^*=\text{argmin}_
{\fontsize{10}{12}\selectfont{\bx_{gray}}}
U_1(\bx)$
and
$\bx_2^*=\text{argmin}_
{\fontsize{10}{12}\selectfont{\bx_{gray}}}
U_2(\bx)$
in grayscale (and color) image restoration problems.
The same idea is used in both cases: it is
to represent of $\bx$ with integer
coordinates $x_i\in \{0,\ldots,L-1\},\ L>2$
by $L-1$ matrices (for general minimization
problem -- by vectors) having $0,1$-coordinates.

\subsection{Efficient minimization of $U_1(\bx)$}

Let for integer $0<l\le L-1$ the indicator functions
${\mathbf 1}_{\ind{(x_i\ge l)}}$ be equal to $1$, if
$x_i\ge l$, and be equal to $0$ otherwise. Then any
grayscale image $\bx$ is represented as the sum
$\bx=\sum_{l=1}^{L-1}\bx(l)$
of binary images-layers $\bx(l)$ with coordinates
$x_i(l)={\mathbf 1}_{\ind{(x_i\ge l)}}$.
Those binary images-layers specify
\emph{the monotone decreasing
sequence} $\bx(1)\ge\bx(2)\ge\ldots\ge\bx(L-1)$, and
vice versa, any sequence of binary images
$\bx(1)\ge\bx(2)\ge\ldots\ge\bx(L-1)$
determines the grayscale image
$\bx=\sum_{l=1}^{L-1}\bx(l)$.
We will use the following obvious
\begin{prop}\label{p:mod}
For any integer $a,b\in \{0,\ldots,L-1\}$ the equality
\begin{equation*}
|a-b|=\sum_{l=1}^{L-1}
\big|{\mathbf 1}_{\ind{(a\ge l)}}-
{\mathbf 1}_{\ind{(b\ge l)}}\big|
\end{equation*}
is valid. Therefore for any grayscale images $\bx,\by$
the function $U_1(\bx)$ can be written in the form
\begin{equation}\label{e:ul}
U_1(\bx)=\sum_{l=1}^{L-1}u_l(x(l)),
\end{equation}
where functions
\begin{equation*}
u_l(x(l))=\la\sum_{i\in S}|y_i(l)-x_i(l)|+
\sum_{(i,j)\in S}\be_{i,j}
\left|x_i(l)-x_j(l)\right|.
\end{equation*}
\end{prop}
Let
\begin{equation}\label{e:lmin}
\check{x}(l)=\text{argmin}_\ind{x_\ind{\text{bin}}}
u_l(x_\ind{\text{bin}}), \quad l\in \{1,\ldots,L-1\}
\end{equation}
denote  binary solutions that minimize the functions
$u_l(x_\ind{\text{bin}})$.
Since $\by(1)\ge\by(2)\ge\ldots\ge\by(L-1)$ the
following sentence is valid.
\begin{prop}\label{p:cx}
There is a monotone decreasing sequence
$\cx_M(1)\ge\cx_M(2)\ge\ldots\ge\cx_M(L-1)$
of solutions
of (\ref{e:lmin}).
\end{prop}
The proof of Proposition~\ref{p:cx}  is similar
to the proof of
Theorem~\ref{t:monot} in Section~\ref{s:MNFC}.
To show existence of solutions $\bx_1^*$ of the
form $\bx_1^*=
\bx_{1,M}^*=\sum_{l=1}^{L-1}\cx_M(l)$ it is enough to
use Propositions \ref{p:mod} and \ref{p:cx}.
Indeed, for any $\bx$
$$
U_1(\bx)=\sum_{l=1}^{L-1}u_l(x(l))\ge
\sum_{l=1}^{L-1}u_l(\cx_M(l))=U_1(\bx_{1,M}^*).
$$
Each binary solution $\cx_M(l)$ can be identified by
the maximum network flow algorithms or by the MNFC for
polynomial number of operations. At worst, it takes
$O(Ln^3)$ operations but quite often the use of MNFC
allows reducing operations up to $O(Ln)$ operations
that are held in a concurrent mode (see
Section~\ref{s:appl}).

\subsection{Efficient minimization of $U_2(\bx)$}

The main idea of the efficient integer minimization of
the function
$U_2(\bx)=\sum_i\la_i(y_i-x_i)^2+
\sum_{i,j}\be_{i,j}(x_i-x_j)^2$
lies in replacement of the variable $\bx$ by the sum of
\emph{unordered Boolean variables} $\bx(l)$, and then
considering another polynomial
$\widetilde{Q}(\bx(1),\ldots,\bx(L-1))
\ge U_2(\sum_{l=1}^{L-1}\bx(l))$
of many Boolean variables
with points of minimum
$\bq(1),\ldots,\bq(L-1)$ such that
$\bx_2^*=\sum_{l=1}^{L-1}\bq(l)$ minimizes $U_2(\bx)$.
The new polynomial $Q$ is chosen so that  it admits the
efficient minimization by the network flow optimization
algorithms.

In more details. Let us represent $\bx$ as the sum of
arbitrary Boolean variables
$\bx=\sum_{l=1}^{L-1}\bx(l)$.
In new variables the polynomial $U_2(\bx)$ will
look as
$$
U_2(\bx)=\sum_i\la_i\lb y_i-\sum_{l=1}^{L-1}x_i(l)\rb^2+
\sum_{i,j}\be_{i,j}\lb\sum_{l=1}^{L-1}(x_i(l)-x_j(l))\rb^2.
$$
Assign for simplicity $b_{i,j}=0$ and use the equality
$x_i^2(l)=x_i(l)$ for Boolean variables $x_i(l)$
to write $U_2(\bx)$ in the form
$$U_2(\bx)=\sum_{i\in S}\la_iy_i^2+P(\bx(1),\ldots,\bx(L-1)),$$
where the polynomial of many Boolean variables
\newpage
\begin{eqnarray}
P(\bx)&=&P(\bx(1),\ldots ,\bx(L-1))\nonumber\\
&=&\sum_{i\in S}\left[\la_i-2\la_iy_i-
(L-2)\sum_{j\in S}(b_{i,j}+b_{j,i})\right]
\sum_{l=1}^{L-1}x_i(l)\nonumber\\
&+&2\sum_{i\in S}\left[\la_i+
\sum_{j\in S}(b_{i,j}+b_{j,i})\right]
\sum_{ 1\le l<m\le L-1}x_i(l)x_i(m)\nonumber\\
&+&\sum_{i,j\in S}b_{i,j}\sum_{l=1}^{L-1}
\bigg( x_i(l)-x_j(l)\bigg)^2\\
&+&\sum_{i,j\in S}b_{i,j}
\sum_{ 1\le l<m\le L-1}
\bigg[(x_i(m)-x_j(l))^2+(x_i(l)-x_j(m))^2 \bigg]\nonumber
\end{eqnarray}
Intoduce another polynomial of many Boolean variables
\begin{eqnarray}
Q(\bx(1),&\ldots& ,\bx(L-1))\nonumber\\
&=&\sum_{i\in S}\left[\la_i-2\la_iy_i-
(L-2)\sum_{j\in S}(b_{i,j}+b_{j,i})\right]
\sum_{l=1}^{L-1}x_i(l)\nonumber\\
&+&2\sum_{i\in S}\left[\la_i+
\sum_{j\in S}(b_{i,j}+b_{j,i})\right]
\sum_{ 1\le l<m\le L-1}x_i(m)\nonumber\\
&+&\sum_{i,j\in S}b_{i,j}\sum_{l=1}^{L-1}
\bigg( x_i(l)-x_j(l)\bigg)^2\\
&+&\sum_{i,j\in S}b_{i,j}
\sum_{ 1\le l<m\le L-1}
\bigg[(x_i(m)-x_j(l))^2+(x_i(l)-x_j(m))^2 \bigg]\nonumber
\end{eqnarray}
such that $Q(\bx(1),\ldots,\bx(L-1))\ge P(\bx)$ and
which differs from $P(\bx)$ by  the term
$\sum_{ 1\le l<m\le L-1}x_i(m)$ in the second summand.

Denote by
\begin{equation*}
(\bq^*(1),\bq^*(2),\ldots,\bq^*(L-1))=
\text{argmin}_{x(1),x(2),\ldots,x(L-1)}
Q(\bx(1),\ldots,\bx(L-1))
\end{equation*}
any collection of Boolean vectors that minimize
$Q(\bx(1),\ldots,\bx(L-1))$.
We will prove that without fail
$\bq^*(1)\ge\bq^*(2)\ge\ldots\ge\bq^*(L-1)$.
This feature will allow expressing solutions of the
initial problem as $\bx_2^*=\sum_{l=1}^{L-1}\bq^*(l)$.
\begin{thm}\label{t:q=x}
Any collection
$(\bq^*(1),\bq^*(2),\ldots,\bq^*(L-1))$
that minimizes $Q$ forms the decreasing sequence.
It identifies the solution of the original problem
by the formula $\bx_2^*=\sum_{l=1}^{L-1}\bq^*(l)$, and
vise versa, each solutions $\bx_2^*$ specifies the
solution $\bq^*(l)=\bx^*(l)$.
\end{thm}
\begin{proof} Let us represent $Q$ as
\begin{multline*}
Q(\bx(1),\ldots,\bx(L-1))=\\
P(\bx)+
\sum_{i\in S}\left[\la_i+
\sum_{j\in S}(b_{i,j}+b_{j,i})\right]
\sum_{ 1\le l<m\le L-1}(1-x_i(l))x_i(m),
\end{multline*}
where $\bx=\sum_{l=1}^{L-1}\bx(l)$. Then
suppose that there exists an unordered collection
$(\bq^*(1),\bq^*(2),
\ldots,\bq^*(L-1))$
minimizing $Q$. For this collection
$$
\sum_{ 1\le l<m\le L-1}
(1- q_i^*(l)) q_i^*(m)>0
$$
and, therefore, for
$\bar{\bx}=\sum_{l=1}^{L-1}\bq^*(l)$
the inequality
$Q\lb\bq^*(1),\ldots,\bq^*(L-1)\rb>
P(\bar{\bx})$
holds
(remind that all
$\la_i,\be_{i,j}>0$). But, evidently,
for ordered Boolean matrices-layers
$\bar{\bx}(1)\ge\ldots\ge\bar{\bx}(L-1)$
with coordinates
$\bar{x}_i(l)=\mathbf{1}_{\ind{(\bar{x}_i\ge l)}}$
we have
\begin{equation*}
Q\lb\bar{\bx}(1),\ldots,\bar{\bx}(L-1)\rb=
P(\bar{\bx}).
\end{equation*}
Hence, unordered sequence $(\bq^*(1),\bq^*(2),
\ldots,\bq^*(L-1))$
can not minimize $Q$.
Since for any ordered collection
$(\bx(1)\ge\ldots\ge\bx(L-1))$
the equality
$Q(\bx(1),\ldots,\bx(L-1))=
P(\bx)$
is valid, actually,
$\bx_2^*=\sum_{l=1}^{L-1}\bq^*(l)$, and
vise versa, each solutions $\bx_2^*$ specifies the
solution $\bq^*(l)=\bx^*(l).$
\end{proof}

It follows from Theorem \ref{t:q=x} that
$P$ and $Q$ have the same ordered set of Boolean
matrices
$(\bx^*(1),\bx^*(2), \ldots,\bx^*(L-1))$
that minimize these polynomials.
But the problem of minimization of $Q$ is known in
discrete optimization
\cite{PR75,Ag83,Ag84}. It is
equivalent to identification of the minimum network
flow cut. To describe the corresponding network let us
re-arrange the polynomial to the form
\begin{eqnarray*}
Q(\bx(1),&\ldots &,\bx(L-1))\nonumber\\
&=&\sum_{i\in S}
\sum_{l=1}^{L-1}\left[(2l-1)\la_i-2\la_iy_i+
(2l-L)\sum_{j\in S}(b_{i,j}+b_{j,i})\right]x_i(l)\\
&+&\sum_{i,j\in S}b_{i,j}\sum_{l=1}^{L-1}
\bigg( x_i(l)-x_j(l)\bigg)^2\\
&+&\sum_{i,j\in S}b_{i,j}
\sum_{ 1\le l<m\le L-1}
\bigg[(x_i(m)-x_j(l))^2+(x_i(l)-x_j(m))^2 \bigg].
\end{eqnarray*}
Then denote for brevity
$$
d_{i,l}=\left[(2l-1)\la_i-2\la_iy_i+
(2l-L)\sum_{j\in S}(b_{i,j}+b_{j,i})\right].
$$
and introduce the network $N(s,t,V,A)$ with the source
$s$, the sink $t$ and usual nodes that are labeled
by $i_l\in V$, where the multiindex $i\in S$ and
$l$ shows a layer number (there are $|V|=(L-1)|S|$ usual
nodes in the $N(s,t,V,A)$). For instance, $i_3$ means
the node $i$ at the third layer. The set of arcs with
corresponding capacities is specified as
\begin{equation*}
A=
\begin{cases}
(s,i_l)  &\\
(i_l,t) &\\
(i_l,j_l)&\\
(i_l,j_m)&\\
\end{cases}
\begin{aligned}
-d_{i,l}        &\ \ \        \text{if}&    d_{i,l}<0,\\
d_{i,l}        &\ \ \       \text{if}&     d_{i,l}>0,\\
b_{i,j}+b_{j,i}&\ \ \        \text{if}&    i\neq j,\\
b_{i,j}+b_{j,i} &\ \ \        \text{if}&  i\neq j,\ l\neq m
\end{aligned}
\end{equation*}
for $l\in\{1,\ldots,L-1\}$. In general, it is
necessary $O(Ln^3)$ operations to find a solution
$\bx_2^*$, but use the MNFC or other decomposition
methods can reduce the operation number.

\section{Applications}\label{s:appl}
It was mentioned in Introduction the MNFC turned out
very efficient in recovering of images (both binary and
grayscale)
corrupted by random noise. In Fig.~1 the binary
$512\times 512$ image
corrupted by rather strong Bernoullian noise with the
parameter $p=0.3$ is placed. The Ising model with the
energy function
$U_1(\bx_{\text{bin}})=U_2(\bx_{\text{bin}})$
was used for image restoration.
\begin{center}
\fbox{\includegraphics[width=6.5cm,height=6.5cm]
{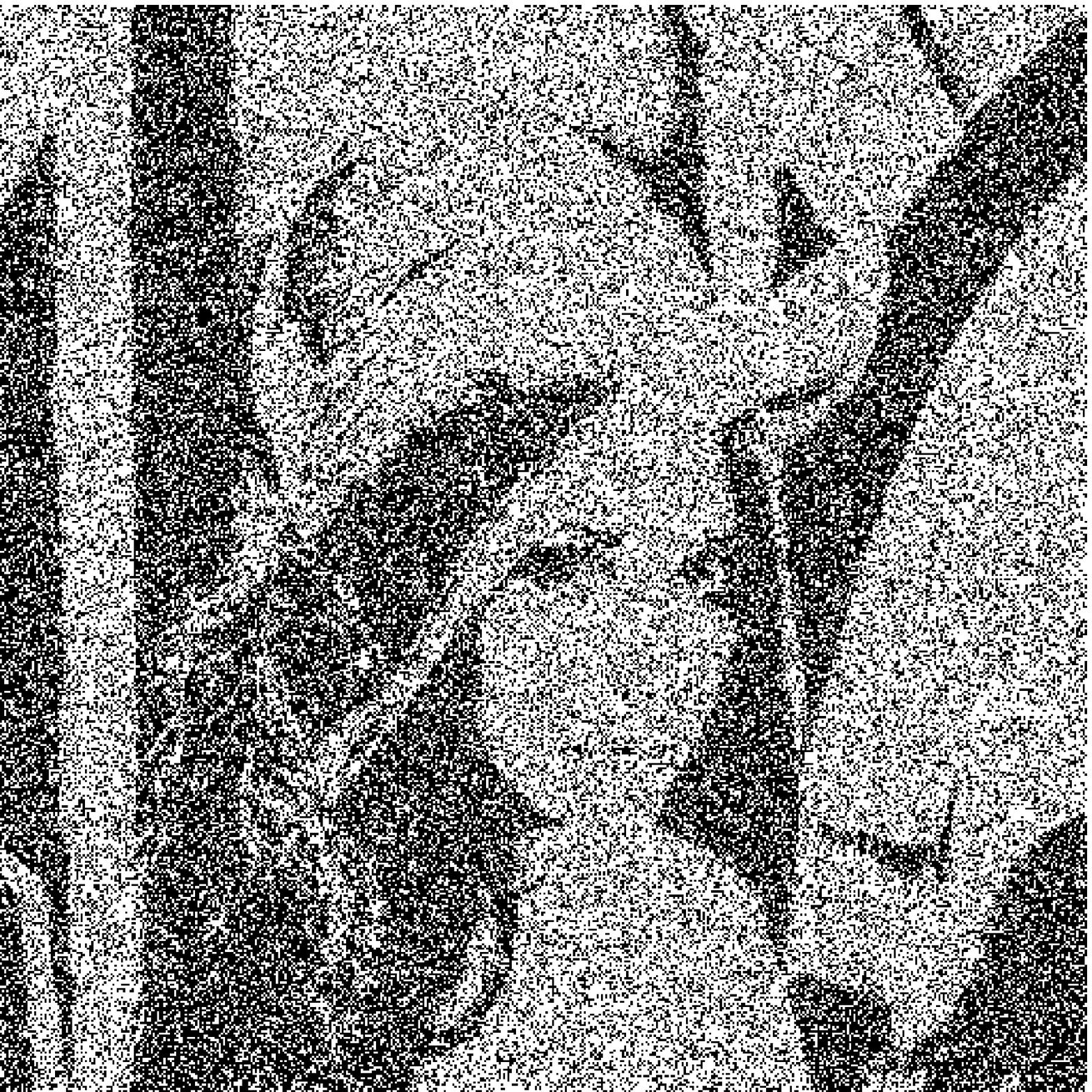}}
\fbox{\includegraphics[width=6.5cm,height=6.5cm]
{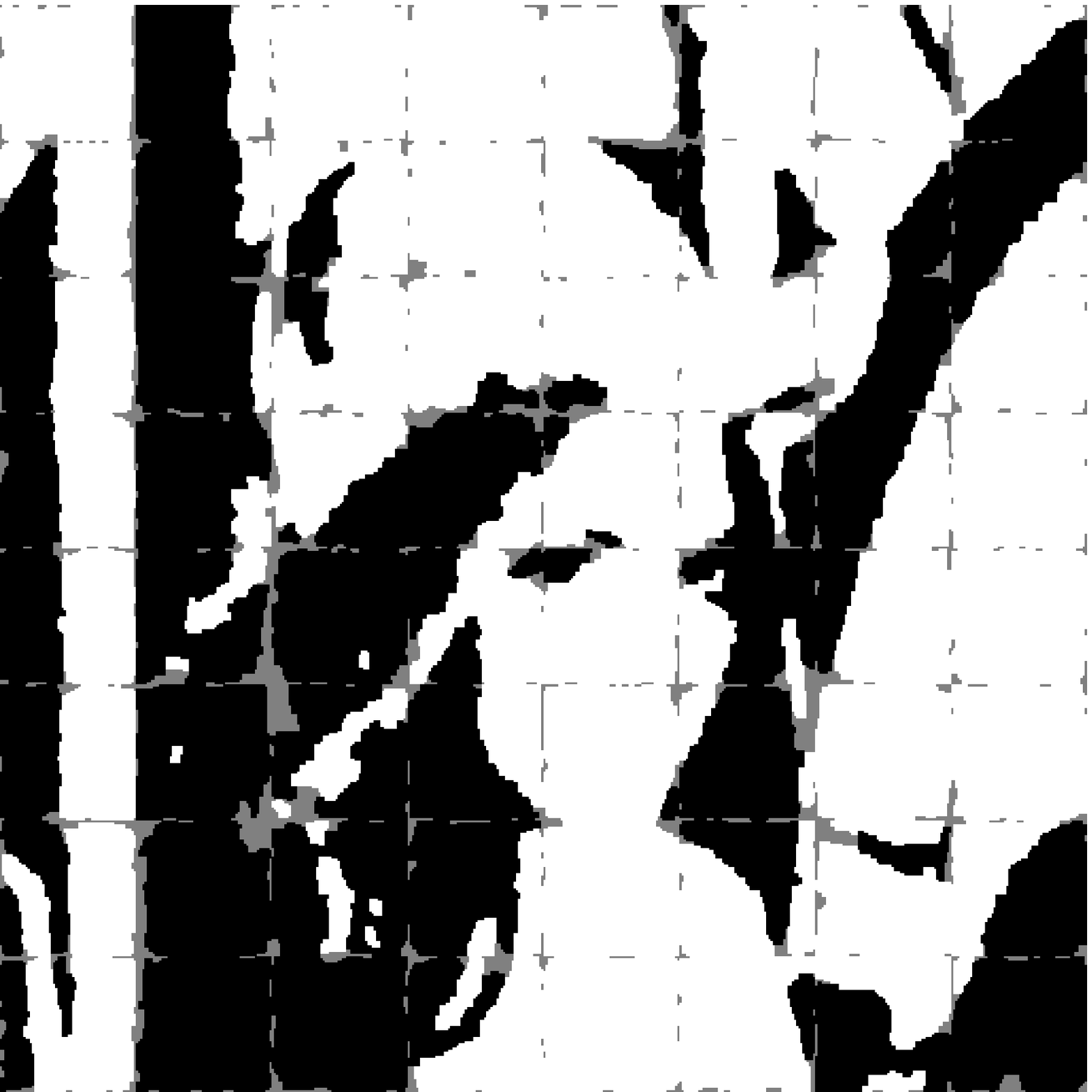}}
\end{center}
\hspace{3cm} Fig.1 \hspace{5.5cm} Fig.2

The result of
implementation of the MNFC at the first level is drawn
in Fig.~2. The pixels that were not estimated at the
first level are depicted in gray color. They form
small sparse sets at boundaries of squares of the partition.
The identification of $\bx^*_{\text{bin}}$ of the
corrupted image takes fractions of a second
while the known maximum
network flow algorithms were not able to evaluate the
estimate at all.

Now consider the application of the proposed method of
the integer minimization for restoration of corrupted
grayscale images. The Ising model with the energy
function $U_1(\bx)$ was used to construct the MAP
estimates. This model was preferred since, first,
the function $U_1(\bx)$ admits more efficient
minimization in comparison with $U_2(\bx)$ and,
second, because of better visual quality of the estimate.
Beside the MAP estimate $\bx_1^*$, the moving average,
the moving median estimates and the gradient estimate
of the $U_2(\bx)$ were computed to compare estimators.

\begin{center}
\fbox{\includegraphics[width=6cm,height=2.4cm]
{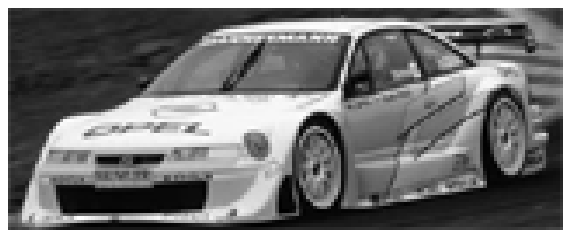}} 
\fbox{\includegraphics[width=6cm,height=2.4cm]
{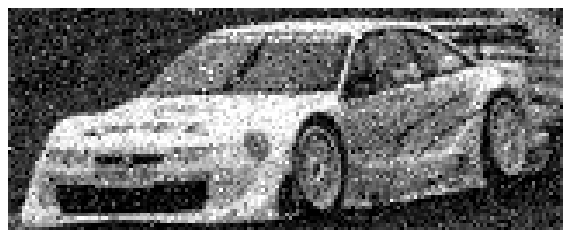}}
\end{center}

\hspace{3.cm}(a)\hspace{6.cm}(b)

\begin{center}
\fbox{\includegraphics[width=6cm,height=2.4cm]
{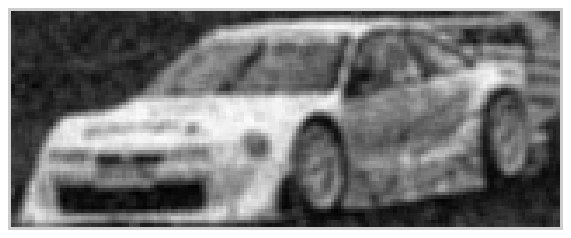}} 
\fbox{\includegraphics[width=6cm,height=2.4cm]
{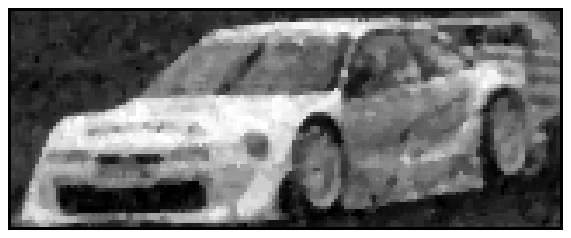}}
\end{center}

\hspace{3.cm}(c)\hspace{6.cm}(d)

\begin{center}
\fbox{\includegraphics[width=6cm,height=2.4cm]
{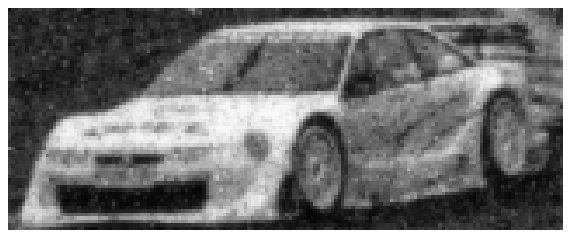}} 
\fbox{\includegraphics[width=6cm,height=2.4cm]
{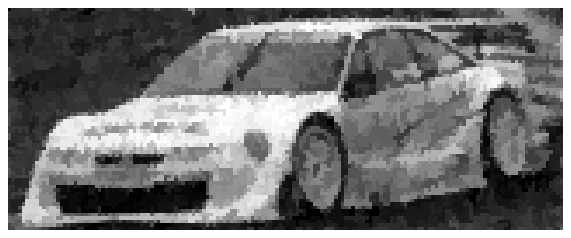}}
\end{center}

\hspace{3cm}(e)\hspace{6cm}(f)
\begin{center}{Fig.~3}\end{center}

The original grayscale image in Fig.~3(a) was corrupted
by the exponential noise Fig.~3(b).
The results of the $3\times 3$ moving average  and
$3\times 3$ moving median filtration are depicted in
Fig.~3(c) and Fig.~3(d) respectively. The continuous
gradient estimate for $\bx_2^*$ is placed in Fig.~3(e).
Fig.~2(f) is the exact $\bx_1^*$ estimate of the image
in Fig.~3(b). Remind that it is determined separately
for every binary layer of a grayscale image. Therefore,
the use of MNFC allows computation of $\bx_1^*$ in a highly
concurrent mode.

\end{document}